\documentclass[12pt,titlepage]{article}
\usepackage{amssymb}
\usepackage[pdftex]{graphicx}
\input epsf
\newtheorem{thm}{Theorem}

\newtheorem{cor}{Corollary}

\newtheorem{prop}{Proposition}

\newtheorem{definition}{Definition}

\newcommand{\text}[1]{\quad\mbox{#1}\quad}
\def\beq{\begin{equation}}\def\eeq{\end{equation}}
\def\beqn{\begin{eqnarray}}\def\eeqn{\end{eqnarray}}

\def\qed{\ifhmode\unskip\nobreak\fi\quad\ifmmode\Box\else$\Box$\fi}

\title{Orientations making $k$-cycles cyclic} 

\author{{\bf Zita Helle}\\
Department of Computer Science and Information Theory\\
Budapest University of Technology and Economics \\
{\tt zita.helle@gmail.com}
\and
{\bf G\'abor Simonyi}\thanks{Research partially supported by the
Hungarian Foundation for Scientific Research Grant (OTKA) Nos.\ K104343 and
K105840.}\\    
MTA Alfr\'ed R\'enyi Institute of Mathematics, \\
{\tt simonyi.gabor@renyi.mta.hu}}

\date{}

\begin{document}
\maketitle

\begin{abstract}

We show that the minimum number of orientations of the edges of the n-vertex
complete graph having the property that every triangle is made cyclic in at
least one of them is $\lceil\log_2(n-1)\rceil$. More generally, we also 
determine the minimum number of orientations of $K_n$ such that at least one
of them orients some specific $k$-cycles cyclically on every $k$-element
subset of the vertex set. The questions answered by these results were
motivated by an analogous problem of Vera T. S\'os concerning triangles and
$3$-edge-colorings.  
\par\noindent
Some variants of the problem are also considered. 
\end{abstract}

\section{Introduction}

How many different sequences of a given length over the alphabet $\{0,1,2\}$ can
be given so, 
that for any three (pairwise distinct) sequences there is a position where
all of them are different, that is exactly one of them is $0$, one is $1$, and
one is $2$? Such a set of sequences we call {\em trifferent}, see
\cite{trif}. Equivalently, what is the minimum length $t$ such that there
exist $M$ distinct ternary sequences forming a trifferent set of
sequences? This is a notoriously hard problem investigated as a special case
of the {\em perfect hashing problem} as defined in \cite{FK, phj}. 

What happens if we want only some specified triples
of our sequences being trifferent? Vera T. S\'os asked such a question in
1991: How many different $3$-colorings of the edges of the complete graph
$K_n$ should be given if for each triangle there is at least one of these
colorings that makes it completely $3$-colored. In \cite{trif} an explicit
construction is given proving that this minimum 
number is at most $\lceil\log_2 n\rceil-1$. The best lower bound we know is
$\lceil\log_3(n-1)\rceil$ given by a simple local consideration:  
the ternary sequences defined by $3$-colorings of the $n-1$ edges connected to
any fixed vertex should all be different. There seems to be no reason to
believe that such a simple lower bound obtained by a local argument would be
sharp.  

In this note we consider an analogous problem involving orientations in place
of colorings. The question here is the minimum number of orientations of the
edges of the complete graph $K_n$ such that each triangle becomes cyclically
oriented in at least one of them. A simple lower bound can be
given here, too, by just looking at orientations locally: If the condition is
satisfied then any pair of the $n-1$ edges attached to a fixed vertex should be
oriented in opposite directions (when viewing them from the common vertex) in
at least one of our orientations. This gives the lower bound
$\lceil\log_2(n-1)\rceil$. 

Somewhat surprisingly, in this case, this simple lower bound is sharp. In
fact, even 
the more general statement of Theorem~\ref{main} below is true.
\medskip
\par\noindent
\begin{definition} Let $n\ge k\ge 3$ and $K_n$ be the complete graph on vertex
  set $[n]=\{0,1,\dots,n-1\}$. A family of 
orientations of the edges of $K_n$ is called {\em increasingly $k$-cycling} if
for 
any $k$-subset $u_1 < u_2 <\dots < u_k$ of the vertex set, there exists an
orientation in our family that makes the cycle $u_1u_2\dots u_ku_1$ cyclically
oriented. 

A family of orientations is called {\em weakly $k$-cycling} if for any
$k$-subset ${w_1,...,w_k}$ 
of $[n]$ there is an orientation making some $k$-cycle consisting of the
vertices $w_1,...,w_k$ (in an arbitrary order) cyclically oriented. 

\par\noindent
The minimum number of orientations in a weakly or increasingly
$k$-cycling family for $K_n$ is denoted by $w(n,k)$ and
$t(n,k)$, respectively. 
\end{definition}
\medskip
\par\noindent
\begin{thm}\label{main}
\begin{equation}\label{exact}
w(n,k)=t(n,k)=\left\lceil\log_2\frac{n-1}{k-2}\right\rceil.
\end{equation}
\end{thm}
\smallskip
\par\noindent
Note that the following corollary immediately follows by setting $k=3$. (It
was proven as a predecessor of the general theorem by the first author in
\cite{HZ}.) Note
also that $t(n,3)=w(n,3)$ is trivial, because all cyclic triangles are ``increasingly
cyclic''. We denote this common value by $t(n)$.  
\smallskip
\par\noindent
\begin{cor} \label{triang} {\rm (\cite{HZ})}
$$t(n)=\left\lceil\log_2(n-1)\right\rceil.$$
\end{cor} 

\medskip\par\noindent
The next section is devoted to the proof of the above theorem, while in
Section~\ref{problems} we discuss some variants of the problem.
\medskip\par\noindent
All logarithms are meant to be of base $2$.

\section{Proof of Theorem~\ref{main}} \label{thmproof}

\subsection{Lower bound}\label{lb}

First we prove the easy fact that the right hand side of (\ref{exact}) is a
lower bound for $w(n,k)$ which obviously bounds $t(n,k)$ from below.

\par\noindent
Assume we have $t':=w(n,k)$ orientations that satisfy the requirements. 
Fix a vertex $u$ and consider the $n-1$ edges adjacent to $u$. Define a binary
sequence of length $t'$ for each of these edges as follows. Let the $i^{\rm
  th}$ bit 
of the sequence belonging to edge $uv$ be $0$ if the edge is oriented from
$u$ to $v$ in the $i^{\rm th}$ orientation. If it is oriented from $v$ to $u$ then
the corresponding bit is $1$. 

We claim that at most $k-2$ edges adjacent to $u$ can get the same $t'$-length
binary
sequence this way. Indeed, if there were $k-1$ edges connecting $u$ to
$v_1,\dots,v_{k-1}$ and all oriented towards $u$ or all oriented away from $u$
then we could not have a cyclic $k$-cycle on vertices $u,
v_1,\dots, v_{k-1}$. This is because that would need a pair of edges $uv_i$
and $uv_j$ 
such that one is oriented towards $u$ and the other away from $u$. Thus we
have at least $\left\lceil\frac{n-1}{k-2}\right\rceil$ different binary sequences, so the length of
these should be at least
$\left\lceil\log\frac{n-1}{k-2}\right\rceil$.

\subsection{Upper bound}\label{ub}

\par\noindent
To prove that the right hand side of (\ref{exact}) is also an upper bound we
give a construction. 

\par\noindent
We may assume that $n=2^r(k-2)+1$ for some positive integer $r$. (Otherwise we
simply make the construction for the first such integer above $n$ and then
delete the superfluous points.) 
  
\par\noindent
Consider the vertices of $K_n$ put around a cycle in the cyclic order
$0,1,\dots n-1,0$. For $a,b\in\{0,1,\dots,n-1\}, a\neq b$ let $d(a,b)$ denote
$b-a-1$ if $b>a$ and $n-2-(a-b-1)=n+b-a-1=n-2-d(b,a)$ if $a>b$. Thus $d(a,b)$
is the number of vertices put on the cycle strictly between $a$ and $b$
when moving from $a$ to $b$ in the
clockwise direction. For each ordered pair $(i,j)$
of vertices we define the value
$f(i,j):=\left\lfloor\frac{d(a,b)}{k-2}\right\rfloor$. (Note that we write
$f(i,j)$ in place of $f((i,j))$.) For each pair of vertices $i\neq j$ consider
the binary form of $f(i,j)$. Since 
$d(a,b)\le n-2<2^r(k-2)$ the number of binary digits of $f(i,j)$ is at most
$r$. 
If it is less, than put as many $0$'s in front of this binary form as are
missing for making it an $r$-digit binary sequence. 
We denote by $f_s(i,j)$ the
$s$th bit in this extended binary form of $f(i,j)$.   
We attach an orientation to each bit, this gives
$r=\log\frac{n-1}{k-2}$ orientations. 
We orient the edge $\{i,j\}$ from $i$ to $j$
in the $\ell$th orientation if $f_{\ell}(i,j)=0$. 
This definition is meaningful since $d(i,j)+d(j,i)=n-2=2^r(k-2)-1$
implying $f(i,j)+f(j,i)=2^r-1$, thus $f_{\ell}(i,j)+f_{\ell}(j,i)=1$ for
all $\ell$, that is $f_{\ell}(i,j)=0$ holds if and only if $f_{\ell}(j,i)=1$. 
\begin{figure}[!htbp]
\centering
\includegraphics[scale=0.14]{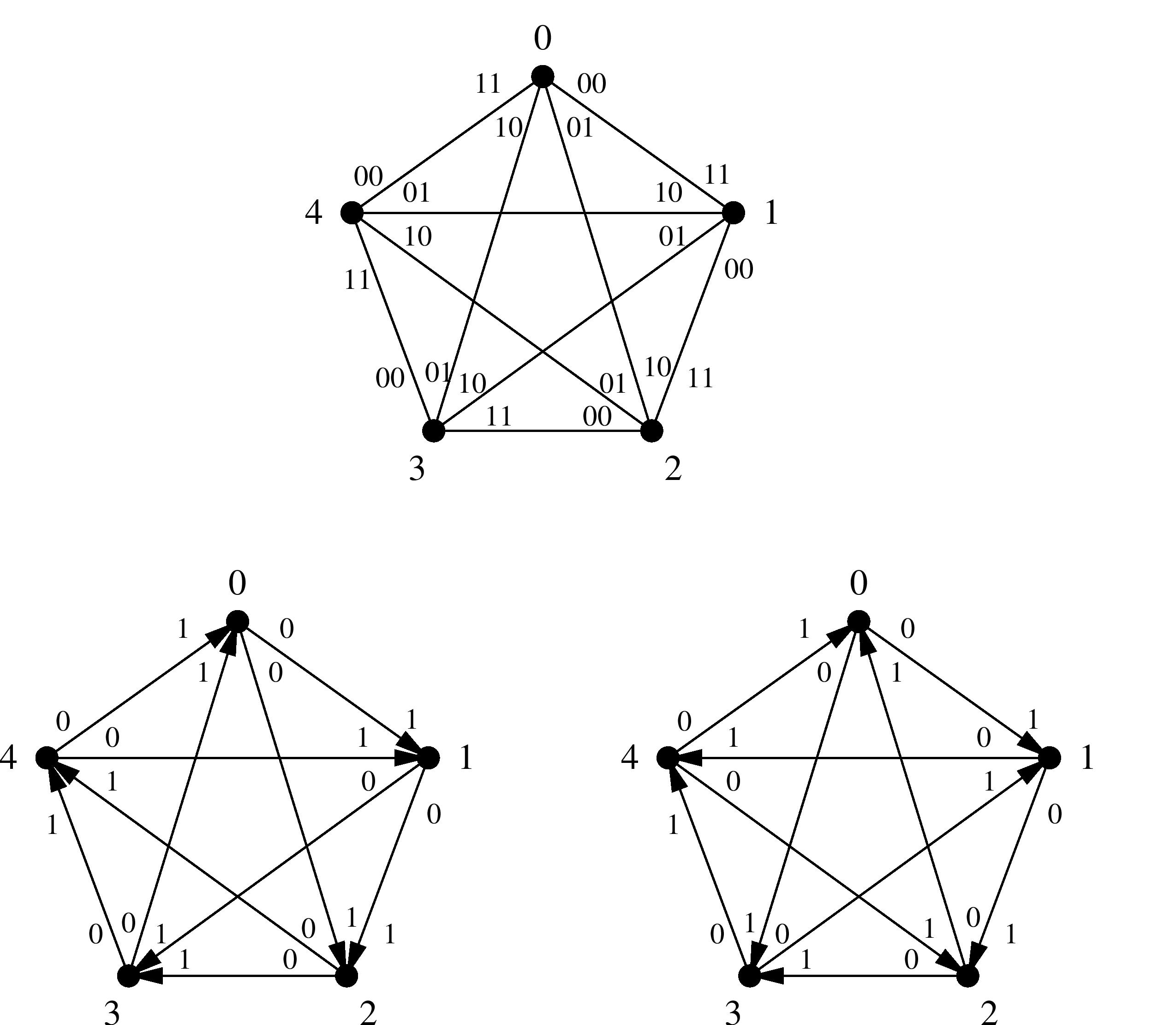}
\caption{The construction for $n=5, k=3$}
\label{fig:figc5}
\end{figure}  
\par\noindent
We have to prove now that this construction works. Let $u_1 < u_2 <\dots < u_k$
be $k$ distinct vertices. The $k$-cycle they form is well oriented in our
$\ell$th orientation if $f_{\ell}(u_i,u_{i+1})=0$ for all $1\le i\le k$, where
the index $k+1$ is meant to be $1$.
To show that such an $\ell\in\{1,\dots,r\}$ exists first note
that $d(u_1,u_2)+d(u_2,u_3)+\dots
d(u_{k-1},u_k)+d(u_k,u_1)=n-k=2^r(k-2)-(k-1).$
Thus we have $\sum_{1\le i\le k}f(u_i,u_{i+1})\le \frac{\sum_{1\le i\le
  k}d(u_i,u_{i+1})}{k-2}<2^r-1$. If there was no $\ell$ for which all
$f_{\ell}(u_i,u_{i+1})=0$, then we would have at least
one $1$-digit at every position making $\frac{\sum_{1\le i<\le
  k}d(u_i,u_{i+1})}{k-2}\ge \sum_{h=0}^{r-1}2^h=2^r-1$
contradicting the previous formula. This contradiction proves the statement.
\hfill$\Box$

\section{Variants and open problems}\label{problems}

There are several ways to generalize or modify the problems treated
above. Here we mention two of them and give some preliminary results.

\subsection{Orienting simplices} 

The following notions often come up in introductory texts on algebraic
topology, in particular when introducing simplicial
homology cf. e.g. \cite{at}. (We mention this only to motivate the following 
definitions, no knowledge of algebraic topology is needed for what follows.)

\smallskip
\par\noindent
\begin{definition}
An ordering of the vertices $v_1,\dots,v_r$ of an $(r-1)$-dimensional simplex
(a simplex with $r$ vertices) 
determines an orientation of the simplex. Two orderings determine the same
orientation if an even permutation of the vertices transforms one to the
other. In particular, there are two possible orientations of a simplex, one
belonging to odd, the other to even permutations of the indices of vertices.
The orientation defined by the permutation $v_{\sigma(1)}\dots v_{\sigma(r)}$
is denoted by $<v_{\sigma(1)}\dots v_{\sigma(r)}>$. The opposite orientation
is denoted by $(-1)<v_{\sigma(1)}\dots v_{\sigma(r)}>$.
\end{definition}

\medskip\par\noindent
{\em Example}:
Let $K$ be a $3$-dimensional simplex with vertices $v_1,\dots,v_4$. The
orderings $v_2v_4v_1v_3$ and $v_3v_2v_1v_4$ belong to the same orientation
(both $2413$ and $3214$ are odd permutations), while $v_4v_3v_2v_1$ determines
the {\em opposite} orientation. So $<v_4v_3v_2v_1>=(-1)<\nolinebreak v_3v_2v_1v_4>$

\medskip
\par\noindent
As usual, the $(d-1)$-dimensional faces of a $d$-dimensional simplex are
called its {\em facets}.

\medskip\par\noindent
\begin{definition}
An orientation $<v_1\dots v_r>$ of the simplex with vertices $v_1\dots v_r$
{\em induces} the orientation $(-1)^i<v_1\dots v_{i-1}v_{i+1}\dots v_r>$ on
its facet $v_1\dots v_{i-1}v_{i+1}\dots v_r$. Two facets $v_1\dots
v_{i-1}v_{i+1}\dots v_r$ and $v_1\dots v_{j-1}v_{j+1}\dots v_r$ are oriented
{\em consistently} if their orientations can be induced by the same
orientation of the original 
simplex. In other words, orientations $<v_1\dots v_{i-1}v_{i+1}\dots v_r>$ and
$<v_1\dots v_{j-1}v_{j+1}\dots v_r>$ are consistent if they induce opposite
orientations on their common face $v_1\dots v_{i-1}v_{i+1}\dots
v_{j-1}v_{j+1}\dots v_r$.
\end{definition}

\medskip\par\noindent
{\em Example}:
Let $K$ be a $3$-dimensional simplex with vertices $v_1,\dots,v_4$. The
orientations $<v_1v_2v_4>$ and $<v_2v_3v_4>$ are consistent, because the first
induces orientation $(-1)<\nolinebreak v_2v_4>$ while the second induces
$(-1)^2<v_2v_4>=<v_2v_4>$ on the face $v_2v_4$ of $K$. Note that the above two
orientations can be the induced orientations of the orientation
$(-1)<v_1v_2v_3v_4>$. 

\medskip\par\noindent
Now the question generalizing the one answered by Corollary~\ref{triang} is
the following: What is the minimum number of orientations of the
($r-2)$-dimensional facets of the
$(r-1)$-dimensional 
simplices defined by $r$-element subsets of $[n]$ if we want
that for every such simplex there is at
least one orientation where all its facets are oriented consistently with each
other? 

\par\noindent
We will use the notation $t_r(n)$ for the number answering the above
question. Note that the function $t(n)$ of Corollary~\ref{triang} is just
$t_3(n)$. 

\medskip\par\noindent
{\em Example}:
For $r=3$ we want to orient the $1$-dimensional (i.e. $2$-element) faces
and for each $2$-dimensional face (i.e. triangle) there is an orientation
where all its facets (the edges) are oriented consistently. The latter means
that at each vertex of the triangle one of the two edges ``comes in'' and the
other ``goes out'', i.e. the corresponding triangle is cyclically
oriented. Corollary~\ref{triang} gives the solution for this case.  

\par\noindent
For $r=4$ we can think about the cyclic orientations of triangles. A $4$-set
($3$-dimesional simplex) is ``satisfied'' if the $4$ triangles it contains get
cyclical orientations so that for all $6$ pairs formed by these triangles, the
common edge of the pair is traversed in the opposite direction by the given
cyclic orientations of these neighboring triangles. 

\medskip\par\noindent
We cannot answer the previous question completely for $r>3$ but give some
estimates. Let us call an orientation of all $(r-1)$-subsets of $[n]$ a {\em
  round}. (That is, it is a function $f: {[n]\choose r-1}\to
\{+,-\}$.) Our problem is to determine the minimum number of
rounds where for each $r$-set there is a round where all its $(r-1)$-subsets
are oriented consistently.

\medskip\par\noindent
\begin{prop} \label{gnr}
$$\lceil\log(n-r+2)\rceil\le t_r(n)\le
\left(\log\frac{2^{r-1}}{2^{r-1}-1}\right)^{-1}\log n+o(\log n).$$ 
\end{prop}

\par\noindent
{\bf Proof.}
The lower bound is a straightforward generalization of the lower bound in
Corollary~\ref{triang}. Fix an $(r-2)$-element subset of $[n]$, we may assume it
is $L:=\{1,\dots,r-2\}$. For an orientation of each $(r-1)$-element subset
containing $L$ we attach a $+$ or $-$ sign according to whether it induces the
orientation $<1,2,\dots,r-2>$ or $(-1)<1,2,\dots,r-2\nolinebreak >$ on
$L$. Since for 
any two $(r-1)$-subsets of $[n]$ containing $L$ we must have a round where one
of them gets a $+$ while the other a $-$ sign, we must have at least $\log
(n-r+2)$ rounds.

\smallskip
\par\noindent
The upper bound is obtained by a simple application of the Lov\'asz
Local Lemma (cf. \cite{ELLL, AS}). We will orient the $(r-1)$-subsets of
$[n]$ randomly in $t$ rounds. Let $A_i$ be the event that the
$(r-1)$-element subsets of the $i^{\rm th}$ $r$-subset of $[n]$ are not oriented
consistently in any of the $t$ rounds. Each  $(r-1)$-subset has two
orientations, so for a fixed $r$-subset there are $2^r$ possible sets of
orientations of its $(r-1)$-element subsets out of which $2$ satisfies our
consistency requirement. So the probability that the requirement is not
satisfied for a fixed $r$-subset in a given round is
$\frac{2^{r-1}-1}{2^{r-1}}$. For $t$ independent rounds we get
$p:=P(A_i)=\left (\frac{2^{r-1}-1}{2^{r-1}}\right )^t$.  
We want to avoid the events  $A_1, \dots, A_{n\choose r}$. A pair of events
$A_j$ and $A_s$ are mutually independent unless the $j^{\rm th}$ and $s^{\rm th}$ $r$-subset of $[n]$
have a common $(r-1)$-subset. So in the dependency graph (formed by the events
$A_ i$ as vertices two being adjacent if they are not mutually independent) the degree of
every vertex is $d:=r(n-r)$. The Local Lemma says (its original form is enough
for our purposes) that if $4pd\le 1$ then all the bad events can be avoided
simultanously. So we need $4\left (\frac{2^{r-1}-1}{2^{r-1}}\right )^tr(n-r)\le
1$. This is satisfied if 
$t\ge\log(4r(n-r))\left(\log\frac{2^{r-1}}{2^{r-1}-1}\right
)^{-1}=\left(\log\frac{2^{r-1}}{2^{r-1}-1}\right)^{-1}\log n+o(\log n).$
\par\noindent
Thus $t_r(n)$ is bounded from above by the smallest $t$ satisfying the
previous inequality and this proves the claimed upper bound.
\hfill$\Box$

\medskip
\par\noindent
Note that Proposition~\ref{gnr} gives $t_4(n)\le 5.20\log n$, while for
$t_3(n)=\lceil\log(n-1)\rceil$ it gives $2.41\log n$ as an upper bound.
\medskip
\par\noindent
{\it Remark 1.} We note that the lower bound in Proposition~\ref{gnr}, though
tight for $r=3$ cannot be tight in general. It is already not so for $r=4$ and
$n=5$. It is not hard to check that we cannot make more than two of the
five $4$-vertex simplices have their facets oriented consistently in one round. Thus at least $3$ rounds are needed, while the lower
bound is $2$. 
\hfill$\Diamond$

\subsection{Orienting every $k$-cycle}

\medskip
\par\noindent
In Section~\ref{thmproof} we considered only some special $k$-cycles of
complete graphs, namely 
those that we called increasing. This property was crucial when proving the
upper bound in Theorem~\ref{main}. Here we look at the problem, where all
$k$-cycles are 
considered. 

\medskip
\par\noindent
Fix a $k\ge 3$. Now our question is the following. What is the minimum number
of orientations of the edges of $K_n$ needed if we want all the $k$-cycles 
being oriented cyclically in at least one of these orientations?  

\medskip
\par\noindent
We do not know the exact answer to this question for $k>3$. Below we show that
this problem is related to the problem of $k$-independent set systems,
cf. \cite{Akind, KS}, and obtain some bounds by applying existing 
results on the latter.  

\par\noindent
\begin{definition}\label{sepsyst}
A family ${\cal F}$ of
subsets of a $t$-element set is called $k$-independent if for every $k$ members
$A_1,\dots, A_k$ of ${\cal F}$, membership of an element in any collection of
these $k$ subsets has no implication to membership in any of the others. In
other words, choosing $B_i$ to be either $A_i$ or its complement, for no such
choice is the intersection $\cap_{i=1}^k B_i$ empty. 
\smallskip
\end{definition}

\medskip
\par\noindent
Kleitman and Spencer \cite{KS} investigated $k$-independent set systems.
Denoting by $f(t,k)$ the maximum number of sets in a $k$-independent family of
subsets of a $t$-element set, they proved that $f(t,k)$ is exponential in
$t$, and in particular, the following upper and lower bounds hold.

\begin{thm} {\rm (Kleitman, Spencer \cite{KS})}\label{KSthm}
For fixed $k\ge 3$ and $t$ sufficiently large we have 
$$2^{d_1t2^{-k}k^{-1}}\le f(t,k)\le 2^{d_2t2^{-k}},$$
where $d_1$ and $d_2$ are absolute constants. 
\end{thm}

If we denote by $T(m,k)$ the minimum size of a set for which an
$m$-element $k$-independent family of its subsets exists, then the above
result of Kleitman and Spencer immediately implies  
\begin{equation}\label{eqkov}
c_12^k\log m\le T(m,k)\le c_2k2^k\log m,
\end{equation}
where $c_1$ and $c_2$ are absolute constants (that can be chosen to be the
reciprocal of $d_2$ and $d_1$, respectively). 

\medskip
\par\noindent
Let $W(n,k)$ denote the minimum size of a family of orientations of the edges of
$K_n$
satisfying that for every odered $k$-element subset $(i_1, i_2, \dots, i_k)$ of
the vertex set there is an orientation in our family which orients the cycle
$i_1\dots i_ki_1$ cyclically (in one or the other direction).

\medskip
\begin{thm}\label{thm:mind}
For $k\ge 3$ fixed and $n$ sufficiently large
$$s_12^{k/2}\log n \le W(n,k)\le s_2k2^k\log n,$$
where $s_1$ and $s_2$ are absolute constants.
\end{thm}

\medskip
\par\noindent
{\bf Proof.}
Let us fix an orientation of the edges of $K_n$ that we will call the {\em
  reference orientation} and denote it by $R$. An arbitrary orientation of
$K_n$ can be encoded by an $n\choose 2$-length sequence of $0$'s and $1$'s: we
attach a $0$ or a $1$ to every edge depending on whether it is oriented the
same way as in $R$ or in the opposite direction.  
Consider a maximum matching
of size $\lfloor n/2\rfloor$. Taking an arbitrary orientation of $\lfloor
k/2\rfloor$ of its edges, it can be extended to a cyclically oriented cycle of
length $k$. Since this cycle (or its reverse) must be present at some
orientation we give, the following must be true. Every possible
$\left\lfloor k/2\right\rfloor$-length $0-1$ sequence
$x=(x_1,\dots,x_{\lfloor k/2\rfloor})$ or its complementary sequence
$\overline x=(1-x_1, \dots, 1-x_{\lfloor k/2\rfloor}))$ must appear in
at least one 
of our $0-1$ sequences describing our orientations, as the
subsequence belonging to the $\lfloor k/2\rfloor$ edges we
looked at. 
This should hold for every choice of $\lfloor k/2\rfloor$ edges
from our $\lfloor n/2\rfloor$ matching edges.
Now double the length of
each of these sequences by concatenating its own
complement to each of them. 
(By the complement of a $0-1$ sequence we mean the sequence
obtained by changing the $0$'s in it to $1$'s and the $1$'s in it to $0$'s.)
Considering the sequences obtained this way as 
characteristic vectors of a set, the size of which is two
times the number of orientations we had, they
must define an $\lfloor n/2\rfloor$-element set of $\lfloor
k/2\rfloor$-independent subsets. This implies $c_12^{k/2}\log (n/2)\le W(n,k)$
by (\ref{eqkov}), and thus the lower bound in the statement. 
\medskip
\par\noindent
For the upper bound consider a $k$-independent
family of $n\choose 2$ subsets of a $T({n\choose 2},k)$-size basic set.
The characteristic vectors of these sets define $T({n\choose 2},k)$ orientations
of the edges of our $K_n$ according to the encoding of orientations described
above. The $k$-independence of our subsets implies that all possible
orientations of any $k$ edges will appear in one of these orientations of
the $K_n$. This holds in particular to the edges of any $k$-cycle, so in one
of the orientations it will be cyclically oriented. Thus $T({n\choose 2},k)$
orientations are sufficient, i.e., $W(n,k)\le T({n\choose 2},k)$. Plugging in
the upper bound of (\ref{eqkov}) we obtain the upper bound in the
statement. This completes the proof. 
\hfill$\Box$

\medskip
\par\noindent
{\it Remark 2.} We note that our problem is also related to
$(i,j)$-separating systems investigated in \cite{FK}.  
For given positive integers $i$ and $j$ a family of bipartitions of a set $X$
forms an $(i,j)$-separating system if for every two disjoint subsets $A$ and
$B$ of $X$ with $|A|=i, |B|=j$ the family contains a bipartition for which $A$
is completely in one, and $B$ is completely in the other partition
class. Fredman and Koml\'os \cite{FK} gave bounds on the minimum number of
bipartitions in an $(i,j)$-separating system of an
$m$-element set. Their results can also be used to prove the lower bound in
Theorem~\ref{thm:mind} and a somewhat weaker upper bound on $W(n,k)$.  
\hfill$\Diamond$

\medskip
\par\noindent
\section{Acknowledgements} 
We are grateful to Benny Sudakov for his inspiring interest. We also thank
G\'abor Tardos for his useful remarks, and in particular for his help in
simplifying the presentation of the proof of Theorem~\ref{main}. A useful
conversation with Kati Friedl is also gratefully acknowledged.

\end{document}